\documentclass[11pt,a4paper]{article}
\usepackage[utf8]{inputenc}
\usepackage[T1]{fontenc}
\usepackage[affil-it]{authblk}
\usepackage{fullpage}
\usepackage[english]{babel}
\usepackage{amsfonts}
\usepackage{latexsym, enumerate}
\usepackage{amsthm, url}
\usepackage{amsopn}
\usepackage{amsmath}
\usepackage[all]{xy}
\usepackage{verbatim}
\usepackage{amssymb}
\usepackage{float}
\usepackage{eucal}
\usepackage[active]{srcltx}
\usepackage{hyperref}
\usepackage{manyfoot}
\usepackage{amscd}
\usepackage[dvips,final]{graphics}
\usepackage[dvips]{graphicx}
\usepackage[usenames,dvipsnames,svgnames,table]{xcolor}
\usepackage{cite}
\usepackage{fancyhdr}

\usepackage[all]{xy}
\usepackage{indentfirst}
\usepackage{fancyhdr}
\usepackage{textcomp}
\usepackage{enumerate}
\usepackage{tikz}
\usetikzlibrary{arrows}
\newtheorem{theo}{Theorem}[section]
\newtheorem{lem}[theo]{Lemma}

\theoremstyle{definition}

\begin{document}
\title{\bf Varieties with Pseudoeffective Canonical Divisor and Their Kodaira Dimension}
\author{Gilberto Bini \thanks{{Universit\`a degli Studi di Palermo, Dipartimento di Matematica e Informatica, Via Archirafi 34, 90123 Palermo (Italia). E-mail: {\tt gilberto.bini@unipa.it}.}}}

\maketitle

\begin{abstract}
Let $X$ be a smooth, projective variety over the field of complex numbers. Here we focus on a conjecture attributed to Shigefumi Mori, which claims that $X$ is uniruled if and only if the Kodaira dimension of $X$ is negative. 
\end{abstract}



\section{Introduction}\label{secintro}

Over the years, the birational classification of smooth, complex projective varieties $X$ has attracted the attention of various algebraic geometers and complex geometers. Their contributions have been so numerous and varied that citing them all would be hardly feasible. Historically, an important chapter was written by Guido Castelnuovo (1865–1952) and by Federigo Enriques (1871–1946) in the first half of the XXth century. Amazingly, they succeeded in classifying smooth projective surfaces $S$ from a birational point of view by means of {\em ad hoc} invariants. In modern terms, the canonical divisor $K_S$ may be either not pseudoeffective or pseudoeffective. In the former case, $S$ is uniruled and the Kodaira dimension of $S$ is negative. In the latter case, the Kodaira dimension of $S$ is nonnegative. Additionally, there exists a {\em good minimal model}, that is, a projective surface $S'$ birational to $S$ such that $K_{S'}$ is semiample, in particular nef. 

To date, the classification of higher dimensional varieties is a wide open problem. In the eighties of the last century, the groundbreaking work of Shigefumi Mori has brought up innovative results, especially in the threefold case: see, for instance, \cite{mor} for an overview. In higher dimensions, there is still considerable work to do in the geography of projective varieties. As in the dimension two case, the canonical divisor may be either not pseudoeffective or pseudoeffective. In the former case, $X$ is uniruled due to the fundamental result in \cite{BDPP}. The latter case is connected to a long-standing conjecture in algebraic geometry, which is attributed to S. Mori and goes back to the beginning of the 1980's, i.e.,  any smooth, complex projective $X$ is uniruled if and only if the Kodaira dimension of $X$ is negative. 

Let us focus momentarily on the claim of Mori's conjecture. Expectedly, if $X$ is uniruled, the Kodaira dimension of $X$ is negative. The difficult implication is the other one, viz., if $kod(X)=-\infty$, the variety $X$ is uniruled or, because of \cite{BDPP}, $K_X$ is not pseudoeffective. The converse of this claim is known as the {\em nonvanishing conjecture}: if $K_X$ is pseudoeffective, the Kodaira dimension of $X$ is nonnegative. Along the lines of the projective surface case, the so-called Abundance Conjecture claims that not only smooth, complex projective varieties but also normal ones with mild singularities admit a good minimal model. This is however beyond the scope of this paper. Here we focus on the Mori conjecture and we prove Theorem \ref{shenzhutaikong} which aims at discussing the nonvanishing conjecture.

Let us briefly summarize the strategy of the proof. Somehow unexpectedly, we manage to find a connection between the Kodaira dimension $kod(X)$ of a smooth, complex projective variety $X$ and the Castelnuovo-Mumford regularity of suitable multiples of $K_X$ via Koszul cohomology groups. In more detail, let $X$ be a $d$-dimensional smooth projective variety over the complex field. Actually, our proof works for $d \geq 3$. Anyhow, the classification of algebraic curves and surfaces is settled, as recalled before. We also assume $K_X$ is pseudoeffective. Pick an ample, globally generated line bundle $A$ on $X$. By abuse of notation, we denote the corresponding Cartier and Weil divisor by $A$. Let $C$ be a generic complete intersection in the numerical class $A^{d-1}$, which is therefore smooth. Denote by ${\mathcal I}_C$ the ideal sheaf of $C$. Since $K_X$ is pseudoeffective, the divisor ${K_X}_{|_C}$ is ample on $C$: see Lemma \ref{taiji}.
Hence, there exists a positive integer $s_A$ such that $s{K_X}_{|_C}$ is effective for every $s \geq s_A$.

Recall that the claim to be proved is the following: if $K_X$ is pseudoeffective, the Kodaira dimension of $X$ is nonnegative. Hence we proceed by contradiction and assume also $kod(X)=-\infty$. For every positive integer $s$ the exact sequence $0 \to {\mathcal I}_C \otimes sK_X \to sK_X\to {\mathcal O}_C \otimes sK_X \to 0$ gives a long exact sequence in cohomology. Moreover, the vanishing of all the plurigenera of $X$ gives a nonzero injective linear map $$H^0(C, {sK_X}_{|_C}) \to H^1(C, {\mathcal I}_C \otimes sK_X)$$ for $s \geq s_A$. Hence $H^1(C, {\mathcal I}_C \otimes sK_X)$ does not vanish as the cohomology group $H^0(C, {sK_X}_{|_C})$ is not zero for $s \geq s_A$.  A spectral sequence argument in Lemma \ref{xinchengmaogong} proves that the Castelnuovo-Mumford regularity $R_0(s):=reg_A\left((1-s)K_X\right)$ of $(1-s)K_X$ with respect to $A$ is greater than or equal to $d+1$. 

As recalled in \eqref{dingtianlinti}, $R_0(s)$ is related to the vanishing of some Koszul cohomology groups. Indeed, $K_{p,q}(X,(1-s)K_X;A)$ vanishes for every integer $p$ and for every integer $q(s) \geq R_0(s)+1$, which implies $q(s)\geq d+2$ because $R_0(s) \geq d+1$, as remarked above. By standard properties of Koszul cohomology,  the following cohomology groups are isomorphic, namely:
$$
K_{p,q}(X,(1-s)K_X;A) \simeq K_{p,1}(X,(1-s)K_X+(q-1)A;A).
$$

Now, write $q(s)=R_0(s)+1 + \delta(s)$ where $\delta(s)$ is a positive integer. Then $(1-s)K_X+(q(s)-1)A$ can be decomposed as the sum of two divisors:
$$
(1-s)K_X+(q(s)-1)A=(\delta(s)A)+ (R_0(s)A+(1-s)K_X)=M_1(s)+M_2(s), \quad s \geq s_A.
$$

Exactly, we prove that $h^0(X,M_i(s))=r_i(s)+1 \geq 2$. Then we apply our technical assumption, which for $A=\mathcal{O}_X$ is a result due to M. Green and R. Lazarsfeld (see \cite{green}, Appendix, Theorem) and prove that for every $s \geq s_A$ there exists an integer $p_{12}(s)$, specifically $p_{12}(s)=r_1(s)+r_2(s)-1$, and an integer $q(s) \geq R_0(s) +1$, such that
\begin{eqnarray*}
\{0\} \neq K_{p_{12}(s), 1}(X, M_1(s)+M_2(s);A) &= & K_{p_{12}(s),1}(X, (1-s)K_X+(q(s)-1)A;A) \\ &= & K_{p_{12}(s),q(s)}(X, (1-s)K_X;A) \\ &= &\{0\}.
\end{eqnarray*}

In the end, if $K_X$ is pseudoeffective, the Kodaira dimension of any smooth, complex projective variety is nonnegative. This is equivalent to saying that if $X$ has negative Kodaira dimension, the canonical divisor $K_X$ is not pseudoeffective, i.e., $X$ is uniruled, as proved in \cite{BDPP}. Vice versa, if $X$ is uniruled, it is well known that the Kodaira dimension is negative. These statements are gathered in Theorem \ref{suntzu}.

\section{Main Results}\label{qi}
 
\subsection{Ample Curves and Castelnuovo-Mumford regularity}
Let $X$ be a $d$-dimensional smooth projective variety over the complex field. As will become clear in what follows, we make the further assumption $d \geq 3$. This is not restrictive due to the classification theory of algebraic curves and algebraic surfaces. We will also assume $K_X$ pseudoeffective. Let $A$ be an ample, globally generated line bundle on $X$.  Let $C$ be a generic complete intersection in the numerical class curve $A^{d-1}$, which is therefore smooth. As remarked in \cite{ottem}, Proposition 4.5., $C$ is an {\em ample curve}.

\begin{lem}
\label{taiji}
Following the notation adopted above, if the numerical dimension $\nu(X)$ of $X$ is different from $0$, the divisor ${K_X}_{|_C}$ is ample on $C$.
\end{lem}
\begin{proof}
Besides being ample, $C$ is a nef subvariety because it is a locally complete intersection (cfr. Proposition 3.6 in \cite{ottem}) and the normal bundle ${\mathcal N}_{C/X}$ is nef, as recalled in \cite{lau}. By Theorem 1.3. in \cite{lau}, the restriction of the pseudoeffective divisor $K_X$ is pseudoeffective. Hence the degree of ${K_X}_{|_C}$ is nonnegative because being pseudoeffective on a curve means being nef. We prove that this degree is in fact positive for every $A$. Assume by contradiction there exists an ample divisor $A$ such that $\deg({K_X}_{|_C})=K_X \cdot A^{d-1}=0$. By the divisorial decomposition of the pseudoeffective divisor $K_X$, we have $(P(K_X)+N(K_X))\cdot A^{d-1}=0$, which implies $P(K_X)\cdot A^{d-1}=0$, as $N(K_X)$ is effective and $P(K_X) \cdot A^{d-1} \geq 0$. Since $C$ is an intersection of ample divisors, this implies that $P(K_X)$ is numerically equivalent to $0$. Hence $\nu(X) = 0$, which is a contradiction, as $\nu(X) \neq  0$ under our assumptions.
\end{proof}

Denote by ${\mathcal F}$ a coherent sheaf on $X$. We recall that the Castelnuovo-Mumford regularity $reg_A({\mathcal F})$ of ${\mathcal F}$ with respect to $A$ (CM-regularity) is defined as
$$
reg_A({\mathcal F})=\min_{m \in {\mathbb Z}} \{h^i(X, {\mathcal F} \otimes A^{m-i})=0, \,\, \forall i >0 \}
.$$ 

The CM-regularity can also be read in terms of the Koszul cohomology groups $K_{p,q}(X, {\mathcal F};A)$, namely:
\begin{equation}
reg_A({\mathcal F})=\min_{m \in {\mathbb Z}} \left \{K_{p,q}(X,{\mathcal F};A)= 0, \,\, \forall p, \,\,\,  \forall q \geq m+1 \right \}.
\label{dingtianlinti}
\end{equation}

\begin{lem}
\label{xinchengmaogong}
Assume $kod(X)= -\infty$. If the CM-regularity of $(1-s)K_X$ is smaller than or equal to $d$ for a positive integer $s$, then $H^1(X, {\mathcal I}_{C} \otimes sK_X)$ vanishes.
\end{lem}
\begin{proof}
The group $H^1(X, {\mathcal I}_{C} \otimes sK_X)$ can be computed in terms of hyper-cohomology. More specifically, the Koszul resolution of the ideal sheaf ${\mathcal I}_{C}$ has $p$-th term given by 
$$
R^{-p}=\bigwedge^{p+1}\underset{d-1}{\underbrace{{\mathcal O}_X(-A) \oplus \cdots \oplus {\mathcal O}_X(-A)}}, \quad 0 \leq p \leq d-2.
$$

Notice that the resolution above exists for $0 \leq p \leq d-2$, whence the assumption $d \geq 3$. There exists a spectral sequence abutting to $H^1(X, {\mathcal I}_{C} \otimes s K_X)$, the $E_1^{p,q}$ term of which is given by the $(p,q)$ term of the double complex $H^q(X, R^{-p} \otimes s K_X)$ for $0 \leq q \leq d$ and $0 \leq p \leq d-2$. Moreover, by Serre's duality we have $E_1^{p,q} = \bigoplus ^{N_p}H^{d-q}(X, (1-s)K_X+(p+1)A)$ where $N_p=\binom{d-1}{p+1}$. Since the spectral sequence degenerates at the $E_2$-page, the vanishing of  $H^1(X, {\mathcal I}_{C} \otimes s K_X)$ originates from the vanishing of $E_1^{p,q}$ for $q=1+p$. This yields $H^i(X, (1-s)K_X+(d-i)A)=0$ for $1 \leq i \leq d-1$. Since the Kodaira dimension is negative, we also have the vanishing for $i=d$. 
\end{proof}

\subsection{Nonnegative Kodaira Dimension}

\begin{theo}
\label{shenzhutaikong}
Let $X$ be a $d$-dimensional smooth, projective variety over the complex numbers for $d \geq 3$. 
Assume that for all divisors $A, M_1, M_2$ on $X$ such that $h^0(X, M_i)=r_i+1, r_i \ge 1, i=1,2,$ we have $K_{r_1+r_2-1,1}(X, M_1 + M_2;A) \neq \{0\}$.
If $K_X$ is pseudoeffective, then the Kodaira dimension $kod(X)$ is nonnegative.
\end{theo}
\begin{proof}
Assume by contradiction $kod(X)=-\infty$. For any positive integer $s$ let us take into account the short exact sequence $0 \to {\mathcal I}_C \otimes sK_X \to sK_X \to {sK_X}_{|_C} \to 0$. By Lemma \ref{taiji}, the divisor ${K_X}_{|_C}$ is ample. As a consequence, there exists a positive integer $s_A$ such that $h^0(C, {sK_X}_{|_C}) \neq 0$ for $s \geq s_A$. Note that here the assumption $K_X$ pseudoeffective is crucial. Since $kod(X)=-\infty$, the linear map $H^0(C, {sK_X}_{|_C}) \to H^1(C, {\mathcal I}_C \otimes sK_X)$ is injective. Hence the cohomology group $H^1(C, {\mathcal I}_C \otimes sK_X)$ does not vanish for any $s \geq s_A$. Therefore, Lemma \ref{xinchengmaogong} implies that $R_0(s)=reg_A\left((1-s)K_X\right) \geq d+1, \quad \forall s \geq s_A.$ The interpretation of $R_0(s)$ in terms of Koszul cohomology groups implies that $K_{p,q}(X, (1-s)K_X;A)=0$ for every integer $p$ and every integer $q(s) \geq R_0(s) +1$. By standard properties of Koszul cohomology, the following vanishing holds:
\begin{equation}
\label{laotzu}
K_{p,1}\left(X, (1-s)K_X+(q(s)-1)A;A\right)=0, \quad \forall p,\quad \forall q(s) \geq R_0(s)+1 \geq d+2, \quad s \geq s_A.
\end{equation}

Notice that $q(s)=R_0(s)+1+\delta(s)$ where $\delta(s)$ is a positive integer, as $s$ varies in the set of positive integers greater than or equal to $s_A$. Indeed, the smallest possible value for $q(s)$ occurs for  $\delta(s)=1$ when $R_0(s)=d+1$. The expression of $q(s)=R_0(s)+1+\delta(s)$ allows us to decompose the divisor $(1-s)K_X+(q(s)-1)A$ as the sum of two divisors:
\begin{equation}
\label{qigong}
(1-s)K_X+(q(s)-1)A=M_1(s)+M_2(s),
\end{equation}
where $M_1(s):=\delta(s)A$ and $M_2(s):= R_0(s)A+(1-s)K_X$. 

\medskip \noindent Note that $h^0(X,M_i(s)) \geq 1$ because $M_i(s)$ are globally generated, in particular effective. Indeed, $M_1(s)$ is so by definition and because $\delta(s) \geq 1$. The line bundle $M_2(s)$ is globally generated because $(1-s)K_X$ is $R_0(s)$-regular: see, for instance, Theorem 1.8.5. (i) in \cite{laz}, volume 1, p. 100. 

\medskip \noindent Moreover, $h^0(X, M_i(s)) \geq 2$. If instead $H^0(X, M_i(s)) = <\sigma_i(s)>$ for a nonzero global section $\sigma_i(s)$, the morphism $0 \to {\mathcal O}_X \to M_i(s) $ defined by $\sigma_i(s)$ is an isomorphism, as $M_i(s)$ is globally generated. On the other hand, $M_1(s)=\delta(s)A$ is not trivial because $\delta(s) \geq 1$, nor is $M_2(s)$; else, $R_0(s)A$ would be equal to $(s-1)K_X$. This would imply $A$ torsion (if $s=1$) or $(s-1)K_X$ effective (if $s \neq 1$), which is impossible under the assumption $kod(X)=-\infty$. 

Now, set $r_i(s)+1=h^0(X, M_i(s))$ for $i=1,2$. As proved above, $r_i(s) \geq 1$; hence we can apply 
our technical assumption inspired by Theorem in \cite{green}, Appendix, and conclude by \eqref{qigong} that $$K_{r_1(s)+r_2(s)-1,1}(X, (1-s)K_X+(q(s)-1)A;A) \neq \{0\}, \quad  \forall s \geq s_A.$$ This however contradicts the vanishing in \eqref{laotzu} for the value of $p$ given by $p_{12}(s)=r_1(s)+r_2(s)-1$. The contradiction originates from the assumption on $X$ that $K_X$ is pseudoeffective and $kod(X)=-\infty$. Thus, any $X$ in the statement of the theorem has nonnegative Kodaira dimension.
\end{proof}

\subsection{Uniruled Varieties}

We collect all the previous results in the following theorem.

\begin{theo}
\label{suntzu}
Let $X$ be a smooth, projective variety over the complex field. Assume that for all divisors $A, M_1, M_2$ on $X$ such that $h^0(X, M_i)=r_i+1, r_i \ge 1, i=1,2,$ we have $K_{r_1+r_2-1,1}(X, M_1 + M_2;A) \neq \{0\}$. Then $X$ is uniruled if and only if $kod(X)=-\infty$.
\end{theo}
\begin{proof}
If $X$ is uniruled, then $kod(X)=-\infty$. Vice versa, Theorem \ref{shenzhutaikong} proves that if $K_X$ is pseudoeffective, $kod(X) \geq 0$. Hence, if $kod(X)=-\infty$, the canonical divisor of $X$ is not pseudoeffective, that is, $X$ is uniruled as proved in \cite{BDPP}. Notice that Theorem \ref{shenzhutaikong} holds for $d$-dimensional varieties for $d \geq 3$. For smaller values of $d$ the claim follows from the classification of algebraic curves and surfaces.

\end{proof}

\end{document}